\documentclass[12pt]{article}
\usepackage{amsmath,amssymb,amsthm,fullpage,hyperref}

\newtheorem{theorem}{Theorem}
\newtheorem{lemma}[theorem]{Lemma}
\newtheorem{corollary}[theorem]{Corollary}
\newtheorem{proposition}[theorem]{Proposition}
\newcommand{\rstir}[3]{\left[ \begin{matrix} #1 \\ #2 \end{matrix}\right]_{#3}}

\AtBeginDocument{}

\author{ \Large{ Matt Davis } \\ \\ Department of Mathematics and Computer Science\\ Muskingum University,  163 Stormont St., New Concord, OH, 43762 \\ mattd@muskingum.edu }
\title{Quadrant Marked Mesh Patterns and the $r$-Stirling Numbers}
\date{}

\begin{document}

\maketitle

\abstract{\textit{Marked mesh patterns} are a very general type of permutation pattern. We examine a particular marked mesh pattern originally defined by Kitaev and Remmel, and show that its generating function is described by the $r$-Stirling numbers. We examine some ramifications of various properties of the $r$-Stirling numbers for this generating function, and find (seemingly new) formulas for the $r$-Stirling numbers in terms of the classical Stirling numbers and harmonic numbers. We also answer some questions posed by Kitaev and Remmel and show a connection to another mesh pattern introduced by Kitaev and Liese.}

\section{Introduction}

The notion of a marked mesh pattern in a permutation is a generalization of classical permutation patterns, and is in fact a common generalization of a number of different variations of permutation patterns that have been of recent interest. Br\"{a}nd\'{e}n and Claesson introduced mesh patterns \cite{BC}, and \'{U}lfarsson  developed marked mesh patterns \cite{U}. We will give a somewhat loose definition of general marked mesh patterns here before turning to the specific examples studied in this paper.

Let $\sigma$ be a permutation written in one-line notation: $\sigma = \sigma_{1}\sigma_{2}\cdots \sigma_{n}$. The graph of $\sigma$ is the set of points of the form $(i, \sigma_{i})$ for $i$ from $1$ to $n$. A \textit{marked mesh pattern} of length $k$ is the graph of a permutation $\sigma$ in $S_{k}$, drawn on a grid, with some regions of the grid demarcated and labeled with a symbol `$=a$', `$\leq a$', or `$\geq a$' for some non-negative integer $a$. This diagram is intended to be a schematic drawing of the graph of a permutation $\sigma$ in $S_{n}$ (with $n>k$), where the drawn points are actual points of the graph of $\sigma$, but the spaces between may have been compressed. Unmarked regions of the pattern place no restriction on $\sigma$, but the marked regions of the graph of $\sigma$ must contain a number of points satisfying the expression that marks the region. (Since the mark ``$=0$'' is the most commonly used, a shaded region with no other mark is assumed to be marked with ``$=0$''.)

The main pattern of interest for this paper is shown in Figure \ref{fig:patterndef}. We denote this pattern as $MMP^{k}$. We precisely define this pattern by saying that $\sigma_{i}$ matches the pattern $MMP^{k}$ in $\sigma$ if there are at least $k-1$ indices $j > i$ with $\sigma_{j} > \sigma_{i}$ and $j < m$, where $\sigma_{m} = n$; and no indices $j < i$ with $\sigma_{j} > \sigma_{i}$. That is, when checking this pattern, we count the points in the graph of $\sigma$ that appear in the rectangle bounded in the lower left by $(i,\sigma_{i})$ and in the upper right by $(m,n)$. For $\sigma_{i}$ to match $MMP^{k}$, there must be at least $k$ points in that rectangle (including the point $(m,n)$ itself, but not $(i, \sigma_{i})$,) and no points above and to the left of $(i, \sigma_{i})$. (This matches the diagram of the pattern because the part of the shaded region along the top of the diagram means that the upper-right point must be $(m,n)$, while the portion of the shaded region on the left gives the restriction on points occurring before $\sigma_{i}$. The middle region specifies the minimum number of points between $\sigma_{i}$ and $n$ in $\sigma$.)

\begin{centering}

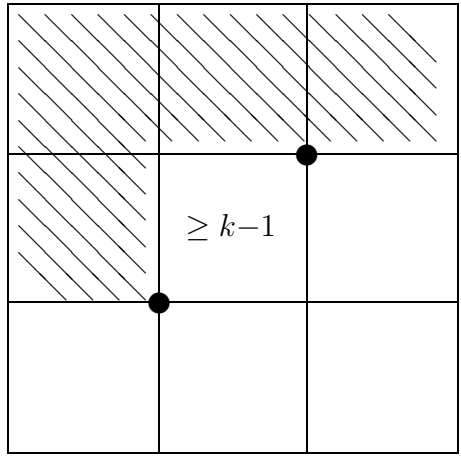
\begin{figure}[ht]
\begin{picture}(400,170)
\put(140,0){\line(1,0){170}}
\put(140,170){\line(0,-1){170}}
\put(140,170){\line(1,0){170}}
\put(310,170){\line(0,-1){170}}
\put(197,170){\line(0,-1){170}}
\put(253,170){\line(0,-1){170}}
\put(140,57){\line(1,0){170}}
\put(140,113){\line(1,0){170}}
\put(197,57){\circle*{8}}
\put(253,113){\circle*{8}}
\put(207,82){$\geq k-$1}
\multiput(144,166)(10,0){12}{\line(1,-1){48}}
\put(264,166){\line(1,-1){38}}
\put(274,166){\line(1,-1){28}}
\put(284,166){\line(1,-1){18}}
\multiput(144,156)(0,-10){6}{\line(1,-1){48}}
\put(144,96){\line(1,-1){38}}
\put(144,86){\line(1,-1){28}}
\put(144,76){\line(1,-1){18}}
\end{picture}

\caption{The Pattern $MMP^{k}$ }

\label{fig:patterndef}
\end{figure}

\end{centering}

The graph of $\sigma$ leads to terminology that will be useful. Draw a horizontal and a vertical line through the point $(i,\sigma_{i})$, which will divide the graph of $\sigma$ into four quadrants. The usual quadrant numbering then lets us more easily describe the position of points relative to $(i, \sigma_{i})$. For example, points above and to the left of $(i, \sigma_{i})$ are said to be in the second quadrant relative to $\sigma_{i}$, or, more succinctly, in $\sigma_{i}$'s second quadrant. Kitaev and Remmel \cite{KR} began a systematic study of patterns for which the restrictions of the pattern can be described in terms of the quadrants relative to $\sigma_{i}$, (including $MMP^{k}$) which they called \textit{quadrant marked mesh patterns}. The main goal of this paper is to describe the generating function for $MMP^{k}$ in terms of $r$-Stirling numbers (allowing us to explain some of the results of Kitaev and Remmel \cite{KR} more fully) and to examine some related combinatorial questions. Specifically, we will find some new recurrences and descriptions for the $r$-Stirling numbers, and note a connection to another marked mesh pattern.

\section{Preliminaries}

For any permutation $\sigma$, we define $mmp^{k}(\sigma)$ to be the number of entries $\sigma_{i}$ that match $MMP^{k}$ in $\sigma$:

\[  mmp^{k}(\sigma) = \left| \{ i \, | \, 1 \leq i \leq n \textrm{ and } i \textrm{ matches } MMP^{k} \textrm{ in } \sigma \} \right|. \]

We will write our permutations in one-line notation. Generally, for $\sigma \in S_{n}$, we use $\sigma_{1}, \sigma_{2}, \ldots, \sigma_{n}$ throughout to denote the entries of $\sigma$. We will also often write $\sigma$ as \begin{equation} \sigma = a_{1,1}a_{1,2} \cdots a_{1,m_{1}}a_{2,1} \cdots a_{2,m_{2}} \cdots a_{p,1} \cdots a_{p,m_{p}}, \label{eq:cycstruct} \end{equation} where $a_{1,1} < a_{2,1} < \cdots < a_{p,1}$, and $a_{q,1} > a_{q,m}$ for any $q$ from 1 to $p$ and $m$ from 1 to $m_{q}$. We will refer to any of the substrings $a_{q,1} a_{q,2} \cdots a_{q,m_{q}}$ as a \textit{pseudocycle} of $\sigma$. The entries $a_{q,1}$ are called \textit{left-to-right maxima} of $\sigma$, since if $\sigma$ is read left-to-right, each entry $a_{q,1}$ is the largest entry read so far. For example, in the permutation $56418732$, the three pseudocycles are the substrings $5, 641$, and $8732$, so we would set $a_{1,1} = 5$, $a_{2,1} = 6$, etc.

We note that the permutations in $S_{n}$ with exactly $k$ pseudocycles are counted by the (unsigned) Stirling numbers of the first kind, denoted $c(n,k)$. These numbers are given by \[ x(x+1)(x+2) \cdots (x+n-1) = \sum_{i=1}^{n} c(n,k)x^{k}.\]

We also define the reduction operator $\textrm{red()}$, which turns a string of $k$ distinct integers into a permutation in $S_{k}$ by setting $\textrm{red}(a_{1}a_{2} \cdots a_{k})$ to be the unique permutation $\sigma \in S_{k}$ resulting from arranging $1$ through $k$ in the same relative order as the numbers $a_{i}$. For example, $\textrm{red}(3625) = 2413$.

\section{Generating Functions}

We now begin our study of the distribution of $mmp^{k}$ on $S_{n}$ with a basic observation relating the number of occurrences of the pattern $MMP^{k}$ in a permutation $\sigma$ to the pseudocycle structure of $\sigma$.

\begin{lemma} \label{lem:basic} If $\sigma \in S_{n}$, then for any $k > 1$, $mmp^{k}(\sigma) \leq \max(0,n-k).$ Moreover, using the notation of \eqref{eq:cycstruct} for $\sigma$, if $mmp^{k}(\sigma) = j$, then the elements that match $MMP^{k}$ in $\sigma$ are exactly $a_{1,1}, a_{2,1}, \ldots, a_{j,1}$. \end{lemma}

Proof. Use the notation of equation \eqref{eq:cycstruct} for $\sigma$ and assume $\sigma$ has $p$ pseudocycles. Then we note that if $\sigma_{i}$ matches $MMP^{k}$ in $\sigma$, then $\sigma_{i}$ is a left-to-right maximum, and thus is one of the $p$ entries $a_{q,1}$ for some $1 \leq q \leq p$. It also implies that there are at least $k$ entries after $\sigma_{i}$, so that $i \leq n-k$. For the second statement, we note that if $a_{i,1}$ matches $MMP^{k}$ in $\sigma$, for $i > 1$, then there are at least $k-1$ elements in the $i+1$st through $p-1$st pseudocycles that are greater than $a_{i,1}$. These same $k-1$ elements occur between $a_{i-1,1}$ and $a_{p,1}$, so that $a_{i-1,1}$ also matches $MMP^{k}$ in $\sigma$. $\square$

The standard generating function for the statistic $mmp^{k}$ is \begin{equation} \label{eq:genfunc} R_{n}^{k}(x) = \sum_{\sigma \in S_{n}} x^{mmp^{k}(\sigma)}.\end{equation} However, a slight alteration to $R_{n}^{k}(x)$ makes its description more straightforward. Let $\sigma = \sigma_{1}\sigma_{2} \cdots \sigma_{n}$. We say that $0$ matches $MMP^{k}$ in $\sigma$ if $n = \sigma_{i}$ for some $i \geq k$. This is really an extension of the definition of $MMP^{k}$, since if we were to prepend a 0 to $\sigma$ and consider it as a permutation $\sigma'$ of $0, \ldots, n$, then $0$ would in fact match $MMP^{k}$ in $\sigma'$ as long as $n$ occurs at the $k$th position or later in $\sigma$. In order to specify when we wish to include this expanded notion of matching $MMP^{k}$, we will define \[ mmp^{k'}(\sigma) = \left| \{ i \, | \, 0 \leq i \leq n \textrm{ and } i \textrm{ matches } MMP^{k} \textrm{ in } \sigma \} \right|. \] Conventionally, we will write $\sigma_{0} = 0$ when we want to consider this leading 0 as part of the permutation. Notice that $mmp^{k'}(\sigma) = 0$ if and only if $n$ occurs at the $k-1$st position or earlier in $\sigma$, and we describe this by saying that $\sigma$ \textit{cannot match} $MMP^{k}$. If $mmp^{k'}(\sigma) = 1$, then we say that $\sigma$ \textit{almost matches} $MMP^{k}$. Note that $n = \sigma_{k}$ is a sufficient but not necessary condition for $mmp^{k'}(\sigma) = 1$. We also point out that if $mmp^{k'}(\sigma) > 0$, then in fact $mmp^{k'}(\sigma) = mmp^{k}(\sigma) + 1$.

As an example that this notion is useful, it can easily be checked that \[ R_{4}^{2}(x) = 17 + 6x + x^{2}.\] However, there are 6 permutations in $S_{4}$ that cannot match $MMP^{2}$, and 11 that almost match $MMP^{2}$. Thus, it seems in some ways more natural to write \[ R_{4}^{2}(x) = 6 + 11 + 6x + x^{2} = c(n,1) + c(n,2) + c(n,3)x + c(n,4)x^{2}.\] More generally, we have the following, which generalizes a result of Kitaev and Remmel \cite[Proposition 3]{KR}.

\begin{theorem} \label{thm:k2} For $s > 0$, the coefficient of $x^{s}$ in $R_{n}^{2}(x)$ is $c(n,s+2)$. The number of permutations that almost match $MMP^{2}$ is $c(n,2)$, while the number of permutations that cannot match $MMP^{2}$ is $c(n,1) = (n-1)!$. \end{theorem}

Proof. First, we claim that if $\sigma$ has $p>2$ pseudocycles, then  $mmp^{2}(\sigma) = p-2$. To see this, write $\sigma$ in cycle notation in the notation of equation $\eqref{eq:cycstruct}$. As in Lemma \ref{lem:basic}, if $\sigma_{i}$ matches $MMP^{2}$ in $\sigma$, then  $\sigma_{i}$ is one of the entries $a_{q,1}$ for some $1 \leq q \leq p$. Conversely, every entry of the form $a_{q,1}$ for $1 \leq q \leq p-2$ matches the pattern $MMP^{2}$ in $\sigma$, since for such $q$, the entries $a_{q+1,1}$ and $a_{p,1} = n$ will lie in the first quadrant relative to $a_{q,1}$. But since $a_{p-1,1}$ is the next-to-last left-to-right maximum, $n$ is the only entry in its first quadrant. Hence, $mmp^{2}(\sigma) = p-2$, as claimed.

For $s > 0$, this implies that the number of permutations $\sigma$ that contribute a term of $x^{s}$ to the sum in Equation \eqref{eq:genfunc} is precisely the number of permutations with $s+2$ pseudocycles, which is counted by $c(n,s+2)$. The count of permutations that cannot match $MMP^{2}$ is immediate from the definition, and the remaining $c(n,2)$ permutations must almost match $MMP^{2}$. $\square$

To more easily keep track of $mmp^{k'}$, we define a new generating function $P_{n}^{k}(x)$ by setting \[ P_{n}^{k}(x) = \sum_{\pi \in S_{n}} x^{mmp^{k'}(\sigma)}. \] For example, \[ P_{4}^{2}(x) = 6 + 11x + 6x^{2} + x^{3}.\] Of course, $R_{n}^{k}(x)$ is easily recoverable from $P_{n}^{k}(x)$ and vice versa, so we will focus our attention on $P_{n}^{k}(x)$. We also define $C_{n,k,j}$ to be the coefficient of $x^{j}$ in $P_{n}^{k}$.

For a fixed $k$, we can generate $P_{n}^{k}(x)$ using the following recursive description:

\begin{theorem} \label{thm:main} For a fixed value of $k > 1$ and any $n \geq k$ and $j > 0$, \[ C_{n,k,j} = (n-1) C_{n-1,k,j} + C_{n-1,k,j-1}.\] That is, \[P_{n}^{k}(x) = (x+n-1) P_{n}^{k}(x).\] Also, \[ P_{k-1}^{k}(x) = (k-1)!. \] \end{theorem}

Proof. First, assume $j > 1$. If $mmp^{k'}(\sigma) = j$, then $\sigma$ must match $MMP^{k}$ at the elements that start each of its first $j-1$ pseudocycles and at 0. We can count such permutations as follows:

If $\sigma_{1} = 1$, then dropping $1$ from $\sigma$ will completely eliminate the first pseudocycle, but will not affect any of the others. Thus $\textrm{red}(\sigma_{2}\sigma_{3} \cdots \sigma_{n})$ will match $MMP^{k}$ at exactly $j-1$ places - the starting point of each of its first $j-2$ pseudocycles and 0. There are exactly $C_{n-1,k,j-1}$ such permutations, and the described correspondence is a bijection, since the inverse map consists of increasing each element of the permutation by 1 and inserting 1 at the start of the string.  If $\sigma_{1} \neq 1$, then $1$ does not match $MMP^{k}$ in $\sigma$, and in fact removing it from $\sigma$ will not affect any of the entries $\sigma_{i}$ that do match $MMP^{k}$ in $\sigma$, with the possible exception of 0. However, since $j >1$, then even after removing $1$ and reducing, $mmp^{k'}(\textrm{red}(\sigma_{1}\sigma_{2} \cdots \hat{1} \cdots \sigma_{n})) = j-1$, so that 0 matches $MMP^{k}$ in $\sigma$ as well. Thus $mmp^{k'}(\textrm{red}(\sigma_{1}\sigma_{2} \cdots \hat{1} \cdots \sigma_{n})) = j$. But, since 1 could be in positions 2 through $n$ in $\sigma$, each permutation $\sigma' \in S_{n-1}$ with $mmp^{k'}(\sigma') = j$ will appear as $\textrm{red}(\sigma_{1}\sigma_{2} \cdots \hat{1} \cdots \sigma_{n})$ for $n-1$ different $\sigma \in S_{n}$. Thus, there are $(n-1) C_{n-1,k,j}$ such $\sigma \in S_{n}$.

If $j=1$, then the theorem claims that the number of permutations in $S_{n}$ that almost match $MMP^{k}$ is equal to $n-1$ times the number of permutations in $S_{n-1}$ that almost match $MMP^{k}$ plus the number of permutations in $S_{n-1}$ that cannot match $MMP^{k}$. This actually can be deduced as a consequence of the result for $j>1$ and the fact that the sum of the coefficients of $P_{n}^{k}(x)$ is $n!$, but we can construct a bijection as well. We assume now that $\sigma$ almost matches $MMP^{k}$. If $n$ appears after the $k$th position of $\sigma$, then $mmp^{k}(\sigma) = 0$ implies that $\sigma_{1} \neq 1$. In this case, or if $n = \sigma_{k}$ and 1 appears after $n$, then the argument from the second case above essentially applies. The only alteration that we must make is to note that under our assumptions, $n$ will still be in the $k$th position or later after 1 is deleted, so that $0$ will still match $MMP^{k}$. This accounts for the term $(n-1)C_{n-1,k,1}$.

For the permutations $\sigma$ where $\sigma_{k} = n$ and 1 appears before $n$, we create a permutation in $S_{n-1}$ that cannot match $MMP^{k}$ as follows - first we swap the positions of $n$ and 1, and then delete 1 and apply $red()$. This is a bijection, since the inverse map consists of increasing each element of $\sigma'$ by 1, inserting 1 in the $k$th position, and swapping the positions of 1 and $n$.

The final statement - the base case for the recursion - is clear from the definition of ``cannot match''. $\square$

A direct consequence of the second and third statements of this theorem is an easy description of $P_{n}^{k}(x)$.

\begin{corollary} \label{cor:genfunc} For $n \geq k-1$, \[ P_{n}^{k}(x) = ((k-1)!) \cdot \prod_{i=k-1}^{n-1} (x+i) .\] \end{corollary}

\section{The \texorpdfstring{$r$}{r}-Stirling Numbers}

The recursive relationship in Theorem \ref{thm:main} is the standard recursion for the (unsigned) Stirling numbers of the first kind. In light of Theorem \ref{thm:k2}, this is not entirely surprising. Of course, for $k > 2$, the initial conditions for $P_{n}^{k}(x)$ are different from those for $c(n,k)$. The resulting coefficients of $P_{n}^{k}(x)$ are a generalization of the Stirling numbers of the first kind defined originally by Mitrinovic \cite{M}, but have appeared in many guises. See Broder \cite{B} and Koutras \cite{Kout} for some examples. We will use the notation $\rstir{n}{m}{r}$ of Broder \cite{B}. These numbers are the $r$-Stirling numbers of the first kind, which are defined as follows. We let $\rstir{n}{m}{r}$ be the number of permutations in $S_{n}$ with exactly $m$ left-to-right-maxima such that $n,n-1, \ldots, n-r+1$ are all left-to-right maxima. (We will typically refer to these as just the $r$-Stirling numbers. While $r$-Stirling numbers of the second kind exist, we will not use them here.) Theorem \ref{thm:main} is enough to show the following result, but our main goal is an explicit combinatorial correspondence connecting $P_{n}(x)$ and the $r$-Stirling numbers.

\begin{theorem} \label{thm:rstir} For $n \geq k$, The coefficient $C_{n,k,j}$ is given by $(k-1)!\rstir{n}{j+k-1}{k-1}$ \end{theorem}

Proof. First, fix $j \geq 1$ and assume $\sigma \in S_{n}$ has exactly $j+k-1$ left-to-right maxima, including $n$, $n-1$, \ldots, $n-k+2$. Say that these left-to-right maxima occur at $\sigma_{i_{1}}, \sigma_{i_{2}}, \ldots, \sigma_{i_{j+k-1}} = n$. We will construct $(k-1)!$ permutations counted by $C_{n,k,j}$, and show that this gives a 1-to-$(k-1)!$ correspondence. We first note that  if $j\geq 1$, $\sigma_{i_{j-1}}$ matches $MMP^{k}$ in $\sigma$, since $\sigma_{i_{j}}, \ldots, \sigma_{i_{j+k-1}} = n$ are all larger than $\sigma_{i_{j-1}}$ and occur after it. If $j=1$ we interpret this conventionally to mean $\sigma_{i_{0}} = 0$ matches $MMP^{k}$ in $\sigma$. Similar logic shows that if $j=0$, then $n$ must occur in position $k=1$ or later in $\sigma$.

Now our argument changes based on whether $\sigma_{i_{j}}$ matches $MMP^{k}$ in $\sigma$. Again, if $j=0$, we will conventionally treat $\sigma_{i_{0}}=0$.

Case 1: First assume $j>0$. If $\sigma_{i_{j}}$ does not match $MMP^{k}$ in $\sigma$, then we claim that any elements $\sigma_{r} \geq \sigma_{i_{j}}$ that occur between between $\sigma_{i_{j-1}}$ and $n$ must be one of $\sigma_{i_{j}}$, $\sigma_{i_{j+1}}, \ldots, \sigma_{i_{j+k-2}}$. If $\sigma_{r}$ were such an ``extra'' element, it could not occur before $\sigma_{i_{j}}$ since $\sigma_{i_{j}}$ is a left-to-right maximum, and it could not occur after $\sigma_{i_{j}}$, or else $\sigma_{i_{j}}$ would match $MMP^{k}$ in $\sigma$. Then if we rearrange the elements $\sigma_{i_{j}}$, $\sigma_{i_{j+1}}, \ldots, \sigma_{i_{j+k-2}}$, but leave all the other entries of $\sigma$ in the same place, then $\sigma_{i_{j-1}}$ will still match $MMP^{k}$ in $\sigma$. Moreover, no matter how these elements are rearranged, none of them will match $MMP^{k}$ in the resulting permutation, since there are only $k-1$ elements between positions $i_{j}$ and $i_{j+k-1}$ that are greater than $\sigma_{i_{j}}$, which is the least of all these elements. Then taking all possible rearrangements of $\sigma_{i_{j}}$ through $\sigma_{i_{j+k-2}}$ gives us $(k-1)!$ permutations $\sigma'$ with $mmp^{k'}(\sigma') = j$.

The analogous case for $j=0$ is the set of permutations with $\sigma_{k-1} = n$. We can rearrange the first $k-1$ entries of $\sigma$ in any order, and all of the resulting permutations cannot match $MMP^{k}$.

Case 2: Now assume $\sigma_{i_{j}}$ does match $MMP^{k}$ in $\sigma$. (Note that we now make no assumptions on $j$, so that case 2 includes the possibility that $j=0$ and $n$ appears at position $k$ or later.) Since the $k-1$ elements $n$, $n-1, \ldots, n-k+2$ must be among the left-to-right maxima in $\sigma$, and $\sigma_{i_{j}}$ is the $k$th largest left-to-right maximum, we must have $\sigma_{i_{j+1}} = n-k+2$, $\sigma_{i_{j+2}} = n-k+3$, etc. This implies that the elements that are larger than $\sigma_{i_{j}}$ that occur between $\sigma_{i_{j}}$ and $n$ are either in the list $\sigma_{i_{j+1}}, \ldots, \sigma_{i_{j+k-1}}$, or occur between some $\sigma_{i_{r}}$ and $\sigma_{i_{r+1}}$. We will create a new permutation $\sigma'$ by rearranging $\sigma$. The entry $\sigma_{i_{j-1}}$ and all entries before it will be left in the same place. However, some of the entries that occur after $\sigma_{i_{j-1}}$ but before $n$ will be moved to the end of the string, according to the following algorithm.

For each element $\sigma_{i_{r}}$, for $r$ from $j+1$ to $j+k-2$, we call an entry that lies between $\sigma_{i_{r}}$ and $\sigma_{i_{r+1}}$ that is larger than $\sigma_{i_{j-1}}$ a ``large follower'' of $\sigma_{i_{r}}$. If $\sigma_{i_{r}}$ has no large followers, then the string $\sigma_{i_{r}} \cdots \sigma_{i_{r+1}-1}$ will not be directly moved when we rearrange $\sigma$. If, however, $\sigma_{i_{r}}$ has at least one large follower, then we will remove the substring starting at $\sigma_{i_{r}}$ and ending at the entry before its final large follower, and move that substring to the end of $\sigma$. If several entries $\sigma_{i_{r}}$ have large followers, then we move the strings starting at each $\sigma_{i_{r}}$ in order based on the last large follower of $\sigma_{i_{r}}$, starting with the smallest one. The permutation that results after we have moved all of these strings will be called $\sigma'$.

For example, let $n=8$, $k=4$, and $j=2$. (Then we are constructing a permutation in $S_{8}$ that matches $MMP^{4}$ exactly once.) If \[ \sigma = 13625748, \textrm{ then } \sigma' = 13548762. \] We have $\sigma_{i_{j}} = 3$, and thus we are in case 2. Both 6 and 7 have large followers, and the algorithm tells us to move the substrings 7 and 62 to the end of the string, in that order, since the remaining large followers are 4 and 5 respectively. As another example, if \[ \sigma = 13647582, \textrm{ then } \sigma' = 13458267,\] since both 6 and 7 are moved to the end of the string (after the 2), in that order.

Now, we show that the  constructed permutation $\sigma'$ satisfies $mmp^{k'}(\sigma') = j$. To do this, it suffices to check that $\sigma_{i_{j}}$ does not match $MMP^{k}$ in $\sigma'$, which will imply that $mmp^{k'}(\sigma') = j$. But, the given algorithm moves all but one large follower from each substring $\sigma_{i_{r}} \cdots \sigma_{i_{r+1}-1}$ after $n$ so that only one element from that substring that is greater than $\sigma_{i_{j-1}}$ remains between $\sigma_{i_{j-1}}$ and $n$. And as noted above, all entries of $\sigma$ greater than $\sigma_{i_{j}}$ that occur between $\sigma_{i_{j}}$ and $n$ are in one of those substrings. Thus $\sigma_{i_{j}}$ matches $MMP^{k-1}$ but not $MMP^{k}$ in $\sigma'$. Then, as in case 1, taking all possible rearrangements of the $k-1$ largest elements that occur after $\sigma_{i_{j-1}}$ and before $n$ gives a total of $(k-1)!$ permutations $\sigma''$ with $mmp^{k'}(\sigma'') = j$, constructed from the original $\sigma$.

Our work so far has given us a way of constructing $(k-1)!$ permutations $\sigma'$ with $mmp^{k'}(\sigma') = j$ out of each permutation $\sigma$ counted by by $\rstir{n}{j+k-1}{k-1}$. Now, we must show that the correspondence we have constructed is exactly one-to-$(k-1)!$ and onto. To do so, we construct a right inverse for the correspondence.

Let $\sigma \in S_{n}$ be a permutation with $mmp^{k'}(\sigma) = j$. We must construct a permutation $\sigma'$ with exactly $j+k-1$ left-to-right maxima, including the set $X = \{n-1, \ldots, n-k+2\}$. Let $s$ be the number of elements from $X$ that occur after $n$ in $\sigma$. By assumption, $\sigma_{i_{j-1}}$, the $j-1$st left-to-right maximum, is the last element that matches $MMP^{k}$ in $\sigma$. Let $a_{1}, a_{2}, \ldots, a_{k-1}$ be the $k-1$ largest elements that occur between $\sigma_{i_{j-1}}$ and $n$. Then we rearrange the $a_{i}$ in such a way that, if we were to switch the positions of the $i$th element (in order as the numbers appear in $\sigma$) of $X$ that appears after $n$ with the $i$th smallest of the $a_{i}$, the elements of $X$ would appear in order.

Finally, for each element $r$ from $n-k+2$ to $n-1$ that occurs after $n$, we remove the substring starting at $r$ and ending just before the next occurrence of an entry from $n-k+2$ to $n-1$. Each of these substrings is then re-inserted just before  so that the substring starting with $n-k+1+i$ is inserted just before $i$.

For example, let $n=7$, $k=4$, and $j=2$. Then given the permutation $\sigma = 1324756$, we have $s=2$, and the $a_{i}$ are the elements 3,2, and 4. We rearrange 2,3, and 4 in $\sigma$ to get $1423756$ since swapping 5 with 2 and 6 with 3 would put 456 in order. Next, we remove the strings 5 and 6 from $1423756$ and re-insert them just before 2 and 3 respectively to get $\sigma' = 1452637$. One can easily check that if we apply the original algorithm to $\sigma'$, the $(k-1)!$ permutations we make from $\sigma'$ are $\sigma$ and the permutations we get by rearranging 2,3, and 4 in $\sigma$.

Because of the way that we arranged the $a_{i}$, each of the elements $n-k+2$ through $n-1$ will be a left-to-right maximum in the resulting permutation. Also, by assumption, there are exactly $j-1$ left-to-right maxima that occur at or before $\sigma_{j-1}$, and $k$ that occur after it (including $n$). Thus, there are exactly $j+k-1$ left-to-right maxima in the resulting permutation. Moreover, this process will reverse the correspondence described above, because of the way (including the order) we chose to move elements with large followers after $n$ in the original algorithm. Thus, the original correspondence is 1-to-$(k-1)!$. $\square$

This description of the coefficients of $P_{n}^{k}(x)$ allows us to answer an unresolved question about a formula for these coefficients. Kitaev and Remmel \cite{KR} observed that a formula for $C_{n,4,2}$ (from entry A001712 in the OEIS \cite{OEIS}) is \[ 6 \sum_{i=2}^{n-3} (-1)^{n+i+1} \binom{i}{2} 3^{i-2} s(n-3,i).\] This formula is, in fact, an example of a formula of Mitrinovic \cite{M}, which we prove a version of here.

\begin{theorem} \label{thm:KRans} For $k \geq 2$, and $0 \leq j \leq n-k+1$, we have \[ \frac{1}{(k-1)!} C_{n,k,j} = \rstir{n}{j+k-1}{k-1} = \sum_{i=j}^{n-k+1} \binom{i}{j} c(n-k+1,i) (k-1)^{i-j}. \] Thus \[ C_{n,k,j} = (k-1)! \sum_{i=j}^{n-k+1} \binom{i}{j} c(n-k+1,i) (k-1)^{i-j}.\] \end{theorem}

Proof. We construct a permutation $\sigma$ counted by $\rstir{n}{j+k-1}{k-1}$. First, we know that the elements $n,n-1, \ldots, n-k+2$ must be left-to-right maxima, and so must occur in order in $\sigma$. We choose a way to insert the entries $1$ through $n-k+1$ to complete $\sigma$. For any $i$ from $j$ to $n-k+1$, choose one of the $c(n-k+1,i)$ permutations of $1, \ldots, n-k+1$ (in one-line notation) with exactly $i$ pseudocycles. We choose $j$ of these pseudocycles to be inserted before $n-k+2$, and insert them in increasing order of their first entry, so that those first entries are left-to-right maxima of $\sigma$. The other $i-j$ pseudocycles will then be inserted somewhere after $n-k+2$: between $n-k+2$ and $n-k+3$, between $n-k+3$ and $n-k+4$, etc., up to between $n-1$ and $n$, or after $n$. That gives $k-1$ possible places that a given pseudocycle could be inserted. If multiple pseudocycles are inserted in the same place, we arrange the pseudocycles in increasing order of their first entries. There are exactly $\binom{i}{j} c(n-k+1,i) (k-1)^{i-j}$ ways to construct $\sigma$ in this way, and every such $\sigma$ can be constructed uniquely in this way. The formula for $C_{n,k,j}$ follows immediately from the first part. $\square$.

Example. Let $n=7$, $k=4$, and $j=2$. Then the theorem gives us \[ 119 = \rstir{7}{5}{3} = 6 \cdot 1 \cdot 9 + 3 \cdot 6 \cdot 3 + 1 \cdot 11 \cdot 1 .\] The first term counts permutations using the permutation 1234 (consisting of 4 different pseudocycles) inserted around 5,6, and 7. The second term includes permutations where 1,2,3, and 4 are arranged in a permutation with three pseudocycles, while the last terms fit 1,2,3,4 into two pseudocycles..

The second part of this theorem is somewhat difficult to see directly. Theoretically, we can apply the correspondence described in Theorem \ref{thm:rstir} to the proof of this theorem to understand it. However, the manner of counting is different for permutations in Case 1 of the proof of Theorem \ref{thm:rstir} as opposed to Case 2. For example, if $n=6$, $k=4$, and $j=2$, the term $\binom{3}{2} \cdot c(3,3) \cdot 3^{1}$ counts the permutations 123564, 124365, 124563, 134256, 134526, 134562, 234156,234516, and 234561. The middle 3 and last 3 of these are easily counted - we choose two numbers out of 1,2,3 to start the permutation and then choose for the remaining one to be placed after 4,5, or 6. However, the first 3 (which correspond to Case 1 in Theorem \ref{thm:rstir}) are counted by choosing 1 and 2 to start the permutation and then choosing which of 3,4, or 5 is placed after 6. Now the factor 3 chooses an element, rather than a position for that element. These two cases make the counting argument much clearer for the $r$-Stirling numbers proper than for the coefficients $C_{n,k,j}$.

\section{A Recurrence For Fixed \texorpdfstring{$n$}{n}}

Now, we examine some more consequences of this description of the coefficients of $P_{n}^{k}(x)$. Broder \cite[Theorem 3]{B} shows that the numbers $\rstir{n}{m}{r}$ satisfy the following recursion in $r$:

\[ \rstir{n}{m}{r} = \frac{1}{r-1} \left( \rstir{n}{m-1}{r-1} - \rstir{n}{m-1}{r}\right) .\]

We prove an equivalent result directly for the coefficients of $P_{n}^{k}(x)$. First, we make an observation.

\begin{lemma} \label{lem:second} Let $\sigma$ be a permutation. For $ j\geq 0$, \[ \textrm{ if } mmp^{k'}(\sigma) =  j+1\textrm{, then } mmp^{(k+1)'}(\sigma) = j+1 \textrm{ or }j.\] \end{lemma}

Proof. The condition of matching $MMP^{k+1}$ is strictly stronger than matching $MMP^{k}$, so $mmp^{(k+1)'}(\sigma) \leq mmp^{k'}(\sigma)$. This implies the result for $j=0$, so we now assume $j \geq 1$. Let $\sigma = a_{1}a_{2} \cdots a_{n}$, and let $a_{i_{1}},a_{i_{2}}, \ldots, a_{i_{j}}$ be the $j$ elements (besides 0) which match $MMP{k}$ in $\sigma$, written so that $i_{1} < i_{2} < \cdots < i_{j}$. Since $a_{i_{j}}$ matches $MMP^{k}$ in $\sigma$, we know that $a_{i_{j-1}} < a_{i_{j}}$, and that there are at least $k$ entries of $\sigma$ larger than $a_{i_{j}}$ that occur after $a_{i_{j}}$ but (weakly) before $n$. Then those $k$ entries and $a_{i_{j}}$ are a set of $k+1$ elements that are larger than $a_{i_{j-1}}$ and occur after it in $\sigma$, but weakly before $n$. Thus $a_{i_{j-1}}$ matches $MMP^{k+1}$ in $\sigma$. Then Lemma \ref{lem:basic} implies that $a_{i_{1}}, \ldots, a_{i_{j-2}}$ match $MMP^{k+1}$, and $mmp^{(k+1)'}(\sigma) \geq j$. $\square$

\begin{theorem} \label{thm:main2} Let $k > 1$ and $j > 0$. Then

\[ (k-1) \cdot \left| \{ \sigma \in S_{n} \, | \, mmp^{k'}(\sigma) = j \textrm{ and } mmp^{(k+1)'}(\sigma) = j-1 \} \right| \]\[ \qquad = \left| \{ \sigma \in S_{n} \, | \, mmp^{k'}(\sigma) = j-1 \textrm{ and } mmp^{(k+1)'}(\sigma) = j-1 \} \right|. \]

As a result, $C_{n,k+1,j-1}$ is given by \[ k \cdot \left| \{ \sigma \in S_{n} \, | \, mmp^{k'}(\sigma) = j \textrm{ and } mmp^{(k+1)'}(\sigma) = j-1 \} \right| .\] \end{theorem}

Proof. We give an explicit $k-1$-to-one correspondence between the two sets in question. Write $\sigma$ in the notation of $\eqref{eq:cycstruct}$, with $a_{1,1} < a_{2,1} < \cdots < a_{p,1}$ the left-to-right maxima of $\sigma$. Moreover, assume that $mmp^{k'}(\sigma) = j$ and $mmp^{(k+1)'}(\sigma) = j-1$. Then, by Lemma $\ref{lem:basic}$, $a_{j-1,1}$ matches $MMP^{k}$ but not $MMP^{k+1}$. (If $j=0$, then $a_{0,1}$ is again considered to be the ``virtual'' entry 0 at the start of $\sigma$.) Then, there are exactly $k-1$ entries of $\sigma$ that are larger than $a_{j-1,1}$ that occur after $a_{j-1,1}$ but strictly before $a_{p,1} = n$. Say that these entries are indexed $a_{i_1}, a_{i_2}, \ldots, a_{i_{k-1}}$, with $a_{j-1,1} < a_{i_1} < a_{i_2} < \cdots < a_{i_{k-1}}$. Note that the $a_{i_q}$ and $a_{j-1,1}$ are the $k$ largest elements that occur in the $j-1$st through $p-1$st pseudocycles of $\sigma$, since the $a_{i_{q}}$ are the only elements larger than $a_{j-1,1}$ that occur in those pseudocycles.

Then, for $q$ from $1$ to $k-1$, define $\pi_{q}$ to be the permutation identical to $\sigma$, except with the positions of $a_{j-1,1}$ and $a_{i_q}$ switched. We claim that $mmp^{k'}(\pi_{q}) = j-1 = mmp^{(k+1)'}(\pi_{q})$. First, we note that the entries $a_{1,1}, a_{2,1}, \ldots, a_{j-2,1}$ still match $MMP^{k+1}$ in $\pi_{q}$, since for any $r$ from $1$ to $j-2$, the set of numbers between $a_{r,1}$ and $a_{p,1}$ in $\pi_{q}$ is the same as in $\sigma$. However, $a_{i_q}$ does not match $MMP^{k}$ in $\pi_{q}$ since the only entries of $\pi_{q}$ that are larger than $a_{i_q}$ that occur after the $j-2$nd pseudocycle but before the last pseudocycle are $a_{i_{q+1}}, \ldots, a_{i_{k-1}}$. However, there are only $k-1-q < k-1$ such elements. Moreover, since $a_{j-1,1}$ began the $j-1$st pseudocycle of $\sigma$, it was greater than all the entries of $\sigma$ before it. Since $a_{i_q} > a_{j,1}$, we know that $a_{i_q}$ is now the first entry of the $j-1$st pseudocycle of $\pi_{q}$. Thus, $mmp^{k'}(\pi_{q}) = j-1 = mmp^{(k+1)'}(\sigma)$.

To finish the proof, we construct the correspondence in the opposite direction. Let $\pi$ be any permutation with $mmp^{k'}(\pi) = j-1$ and $mmp^{(k+1)'}(\pi) = j-1$. Write \[ \pi = a_{1,1}a_{1,2} \cdots a_{1,m_{1}}a_{2,1} \cdots a_{2,m_{2}} \cdots a_{p,1} \cdots a_{p,m_{p}},\] with the usual notation. Then, by assumption and Lemma $\ref{lem:basic}$, $a_{j-2,1}$ matches $MMP^{k+1}$ and $MMP^{k}$ in $\pi$, but $a_{j-1,1}$ matches neither. Thus, there are at least $k$ entries of $\pi$ larger than $a_{j-2,1}$ that occur after $a_{j-2,1}$ but strictly before $a_{p,1} = n$, but at most $k-1$ of those occur after $a_{j-1,1}$. Since none of those entries larger than $a_{j-2,1}$ can occur in the $j-2$nd pseudocycle, there must be exactly $k$ of these - $a_{j-1,1}$ and $k-1$ others that occur after $a_{j-1,1}$ - which we call $a_{i_1}, a_{i_2}, \ldots, a_{i_k}$. Assume that $a_{i_1} < a_{i_2} < \cdots < a_{i_{k}}$. Let $q$ be the index so that $a_{j,1} = a_{i_q}$. Notice that $q \neq 1$, since if it were, then $a_{j,1}$ would match $MMP^{k}$ in $\pi$. Then, let $\sigma$ be the permutation identical to $\pi$, except with the positions of $a_{i_1}$ and $a_{i_q}$ switched. Then by construction, $\pi$ is exactly the permutation $\pi_{q}$ constructed from $\sigma$ as in the previous paragraph.

This shows that the correspondence above is exactly $(k-1)$-to-one, and its image is the entire set $\{ \sigma \in S_{n} \, | \, mmp^{k'}(\sigma) = j-1 \textrm{ and } mmp^{(k+1)'}(\sigma) = j-1 \}$.

The second part of the result then follows from the first part and Lemma \ref{lem:second}. $\square$

This result can be seen in the following table displaying the generating functions $P_{5}^{k}(x)$, beginning with $k=2$ in the first row. The arrows between the polynomials demonstrate the number of permutations in the sets described by the theorem. For example, the arrow from $10x^3$ in the first row to $18x^2$ in the second row depicts the 9 permutations $\sigma$ with $mmp^{2'}(\sigma) = 3$ and $mmp^{3'}(\sigma) = 2$. We can determine the labels on the arrows from right-to-left, since the first diagonal arrow on the right in each row must be labeled with the coefficient of the highest degree term, each diagonal arrow determines the vertical arrow pointing to the same term (by the theorem), and the total labels on the two arrows coming from any given term must add up to the coefficient of that term.

\begin{figure}
\begin{centering}
$\begin{array}{ccccccccc}
24 & \, & + 50x & \, & + 35x^2 & \,& + 10x^3 & \, & + x^4  \\
\downarrow_{24} & \swarrow_{24} & \downarrow_{26} & \swarrow_{26} & \downarrow_{9} & \swarrow_{9} & \downarrow_{1} & \swarrow_{1} & \,  \\
48 & \, & + 52x & \, & + 18x^2 & \, & +2x^3 & \, & \, \\
 \downarrow_{48} & \swarrow_{24} & \downarrow_{28} & \swarrow_{14} & \downarrow_{4} & \swarrow_{2} & \, & \, & \, \\
  72 & \, & +42x & \, & + 6x^2 & \, & \, & \, & \, \\ \downarrow_{72} & \swarrow_{24} & \downarrow_{18} & \swarrow_{6} & \, & \, & \, & \, \\ 96 & \, & +24x & \, & \, & \, & \, & \, & \, \\ \downarrow_{96} & \swarrow_{24} & \, & \, & \, & \, & \, & \, & \, \\ 120 & \, & \, & \, & \, & \, & \, & \,
\end{array}
$
\caption{A Recursive Calculation of $P_{5}^{k}(x)$ for $k \leq 5$. \label{fig:extable}}

\end{centering}
\end{figure}

We define the notation \[ m_{n,k,i,j} = \left| \{ \sigma \in S_{n} \, | \, mmp^{k}(\sigma) = i \textrm{ and } mmp^{k+1}(\sigma) = j \} \right| \] for the marks on the arrows in diagrams like this one.

As a result, we can easily deduce the following recurrence relation satisfied by the coefficients of $P_{n}^{k}(x)$. It essentially describes the process of determining one row of the above table from the previous row.

\begin{corollary}  For $k \geq 2$ and $0 \leq j \leq n-k$, \[ C_{n,k+1,j} = \sum_{i=j+1}^{n-k+1} k(1-k)^{i-j-1} C_{n,k,i} . \]
\end{corollary}

Proof. We prove the corollary inductively. By Lemma \ref{lem:second}, if $mmp^{k'}(\sigma) = n-k+1$, then $mmp^{(k+1)'}(\sigma) = n-k$. Thus the theorem implies that $C_{n,k+1,n-k} = k \cdot C_{n,k,n-k+1}$. That is, the corollary holds for $j = n-k$. Then, assume that the formula in the corollary holds for some $j > 1$. This implies (using the second part of the theorem) that \begin{equation} \label{eq:corIH} \sum_{i = j+1}^{n-k+1} (k)(1-k)^{i-j-1} C_{n,k,i} = C_{n,k+1,j} =  k \cdot m_{n,k,j+1,j}. \end{equation} Then, by Lemma \ref{lem:second} and the theorem, \begin{equation} \label{eq:corsec} C_{n,k+1,j-1} = m_{n,k,j-1,j-1} + m_{n,k,j,j-1} = k \cdot m_{n,k,j,j-1} \end{equation} and \begin{equation} \label{eq:corthird} m_{n,k,j,j-1} = C_{n,k,j} - m_{n,k,j,j} = C_{n,k,j} - (k-1) m_{n,k,j+1,j}.\end{equation} Thus, by \eqref{eq:corIH}, \eqref{eq:corsec}, and \eqref{eq:corthird}, \[ C_{n,k+1,j-1} = k(C_{n,k,j} + (1-k)m_{n,k,j+1,j}) \] \[ = kC_{n,k,j} + (1-k) \cdot \sum_{i = j+1}^{n-k+1} (k)(1-k)^{i-j-1} C_{n,k,i} \] \[ = \sum_{i = j}^{n-k+1} (k)(1-k)^{i-j} C_{n,k,i}. \] That is, the theorem is true for $j-1$, and by induction, it is true for all $j \geq 0$.  $\square$

\begin{corollary} \label{cor:sumform} For $k \geq 2$ and $0 \leq j \leq n-k$, \[ \rstir{n}{j+k}{k} = \sum_{i=j+1}^{n-k+1} (1-k)^{i-j-1} \rstir{n}{i+k-1}{k-1}. \] \end{corollary}

Proof. Simply substituting the formula for $C_{n,k,j}$ given in Theorem \ref{thm:rstir} into the previous corollary and canceling $k!$ from both sides gives this result. $\square$

This corollary is a special case of a theorem of Broder \cite[Theorem 19]{B}, but is included since it is such an easy consequence of Theorem \ref{thm:main2}. However, our main goal for the results of this section (and pictures like Figure \ref{fig:extable}) is to describe a (seemingly new) connection between the $r$-Stirling numbers and the classical Stirling numbers.

\section{The Classical Stirling Numbers}

The recursive description of the coefficients of $P_{n}^{k}(x)$ gives us another way to compute an explicit formula for any coefficient $C_{n,k,j}$. For example, we can compute a formula for the coefficients of $P_{n}^{3}(x)$ as follows. We assume that $n$ is large, and begin by writing the entries $c(n,1)$, $c(n,2)$, etc. in the first row. We can then mark the arrows in the diagram left-to-right - the first vertical arrow (And thus the first diagonal arrow) is marked with $c(n,1)$. Thus the second pair of arrows are both marked with $c(n,2) - c(n,1)$. Continuing in this fashion, we can fill out the entire second row of the diagram.

\begin{figure}[ht]
\begin{centering}

$\begin{array}{cccccc}
c(n,1) & \, & c(n,2)x & \, & c(n,3)x^{2} & \cdots \\
\downarrow_{c(n,1)} & \swarrow & \downarrow_{c(n,2)-c(n,1)} & \swarrow & \downarrow_{c(n,3)-c(n,2)+c(n,1)} & \cdots \\
2c(n,1) & \, & 2(c(n,2) - c(n,1))x & \, & 2(c(n,3)-c(n,2) + c(n,1))x^{2} & \cdots
\end{array}$
\caption{Table to compute $C_{n,3,j}$ in terms of $c(n,k)$.}

\end{centering}
\end{figure}

The alternating sum of Stirling numbers of the first kind that appear in this diagram is in fact plus or minus the sum of signed Stirling numbers of the first kind, whichever makes the total sum positive. Thus, the computations depicted in the diagram prove:

\begin{corollary} \label{k3cor}

For any $n > 1$, \[ P_{n}^{3}(x) = \sum_{j=0}^{n-2} 2\left| \sum_{i=1}^{j+1} s(n,i) \right| x^{j} .\]

\end{corollary}

Kitaev and Remmel pointed out a special case of this theorem \cite[Proposition 4]{KR}. As another example, we find a formula for the number of permutations that almost match $MMP^{k}$. To do so, we fill out only the first two columns of the diagram above. The first vertical and diagonal arrows are determined, and from the first diagonal arrow and the second entry in each row, we can compute the mark on the second vertical arrow. Thus we can compute the entire table one column at a time, since each entry in the $k$th row of the table is $\frac{k}{k-1}$ times the mark on the vertical arrow pointing to it.

\begin{figure}[ht]

\begin{centering}
$\begin{array}{ccccc}

c(n,1) & \, & c(n,2)x  & \, & c(n,3) x^{2}\\
\downarrow_{c(n,1)} & \swarrow & \downarrow_{c(n,2) - c(n,1)} & \swarrow & \downarrow_{c(n,3) - c(n,2) + c(n,1)} \\
2c(n,1) & \, & \left(2c(n,2)-2c(n,1)\right)x  & \, & \left(2c(n,3) - 2c(n,2) + 2c(n,1)\right)x^{2} \\
\downarrow_{2c(n,1)} & \swarrow & \downarrow_{2c(n,2) - 3c(n,1)} & \swarrow & \downarrow_{2c(n,3) - 3c(n,2) + \frac{7}{2}c(n,1)} \\
3c(n,1) & \, & \left(3c(n,2) - \frac{9}{2}c(n,1)\right)x & \, & \left(3c(n,3) - \frac{9}{2}c(n,2) + \frac{21}{4} c(n,1) \right)x^{2} \\
\downarrow_{3c(n,1)} & \swarrow & \downarrow_{3c(n,2) - \frac{11}{2}c(n,1)} & \swarrow & \downarrow_{3c(n,3) - \frac{11}{2}c(n,2) + \frac{170}{24}c(n,1)} \\
4c(n,1) & \, & \left(4c(n,2) - \frac{44}{6}c(n,1)\right)x & \, &  \left( 4c(n,3) - \frac{44}{6}c(n,2) + \frac{340}{36} c(n,1) \right) x^{2} \\
\downarrow_{4c(n,1)} & \swarrow & \downarrow_{4c(n,2) - \frac{50}{6}c(n,1)x} & \swarrow & \downarrow_{4c(n,3) - \frac{50}{6}c(n,2) + \frac{415}{36}c(n,1)} \\
5c(n,1) & \, & \left(5c(n,2) - \frac{250}{24}c(n,1)\right) x & \, & \left( 5c(n,3) - \frac{250}{24}c(n,2) + \frac{2075}{144}c(n,1) \right)x^{2}
\end{array}$

\caption{Table to compute $C_{n,k,1}$ and $C_{n,k,2}$ in terms of $c(n,k)$ \label{fig:seccolform}}

\end{centering}

\end{figure}

We can see from this computation that the number of permutations that almost match $MMP^{k}$ is of the form $(k-1)c(n,2) - Ac(n,1)$ for some positive constant $A$, which can be found recursively, or using boundary conditions to solve for $A = \frac{c(k,2) -(k-2)!}{(k-2)!}$. (The numerator of the coefficient of $c(n,1)$ in this formula is A052881 in the OEIS \cite{OEIS}.)

In light of Theorem \ref{thm:rstir}, this formula expresses $\rstir{n}{r+1}{r}$ as a linear combination of the classical Stirling numbers. But this procedure could be continued, so that we have the following:

\textbf{Fact}: For any $j \geq 0$ and $k \geq 2$, there are constants $a_{m,r,1}, a_{m,r,2}, \ldots, a_{m,r,m-r}$ so that, for all $n > j + k -2$, we have \[ \rstir{n}{m}{r} = \sum_{i=1}^{m-r} a_{m,r,i} c(n,i) .\]

We will describe these coefficients more precisely in Theorem \ref{thm:lincomb} below. Note that the $r$-Stirling numbers were previously known to be linear combinations of classical Stirling numbers of the first kind. Koutras \cite[Equation 1.8]{Kout} (for example) describes $\rstir{n}{m}{r}$ as a linear combination of $c(k,k), c(k+1,k), \ldots, c(n,k)$ rather than $c(n,1),c(n,2), \ldots, c(n,k)$, as we have here. Broder \cite[Equation 43]{B} also describes $\rstir{n}{m}{r}$ as a linear combination of Stirling numbers, but as an integer linear combination of $c(n,m), c(n,m+1), \ldots, c(n,n)$, so that the formulas described in this fact are new. (In fact, Broder's formula can be recovered by drawing a table like Figure \ref{fig:seccolform} that starts from the right-hand side rather than the left.)

\section{Harmonic Sums}

The Stirling numbers (and thus the $r$-Stirling numbers) have a number of connections to harmonic numbers, which we now explore.  We define the \textit{$n$th harmonic number} to be $H_{n}^{(1)} = \sum_{i=1}^{n} \frac{1}{i}$. We recursively define the \textit{iterated harmonic sum} for $j > 1$ by \[ H_{n}^{(j)} = \sum_{i=1}^{n} \frac{H_{i}^{(j-1)}}{i} .\] Conventionally, we will set $H_{n}^{(0)} = 1$ for any $n > 0$.

We will be most interested in the sequences $n H_{n-1}^{(j)}$. For example, $n H_{n-1}^{(1)}$ is the sequence \[ 2,9/2, 44/6, 250/24, \ldots,\] which are exactly the coefficients of $c(n,1)$ that appear in the second column of Figure \ref{fig:seccolform}. We can also check that $n H_{n-1}^{(2)}$ is \[ 2, 21/4, 340/36, \ldots, \] which appear in the third column of Figure \ref{fig:seccolform}. This leads us to the following theorem:

\begin{theorem} \label{thm:lincomb} For $j>0$, and $k>2$, the coefficients of $P_{n}^{k}(x)$ satisfy \[ C_{n,k,j} = \sum_{i=1}^{j+1} (k-1) H_{k-2}^{(j+1-i)} (-1)^{j+1-i} c(n,i).\] \end{theorem}

Proof. We proceed inductively. For $k=3$, the statement is exactly the same as Corollary \ref{k3cor}. Now, we assume the result for a fixed value of $k$, and prove the result for $C_{n,k+1,j}$ by induction on $j$. Again, the result for $j=0$ is obvious since $C_{n,k+1,0} = k(n-1)! = k c(n,1)$, as the theorem claims. Now, we assume the result for $C_{n,k+1,j}$ for some fixed $j$, and prove the result for $j+1$.

By assumption, $C_{n,k+1,j} = \sum_{i=1}^{j+1} k H_{k-1}^{(j+1-i)} (-1)^{j+1-i} c(n,i).$ Then, we imagine the tables as drawn in Figure \ref{fig:seccolform}, focusing on $C_{n,k,j+1}$, $C_{n,k+1,j}$, and $C_{n,k+1,j+1}$. Recall the notation $m_{n,k,j,j'}$ for the mark on the arrow pointing from $C_{n,k,j}$ to $C_{n,k+1,j'}$. By Theorem  \ref{thm:main2}, $m_{n,k,j+1,j} = \sum_{i=1}^{j+1} H_{k-1}^{(j+1-i)} (-1)^{j+1-i} c(n,i).$ Then, we compute: \begin{align*} m_{n,k,j+1,j+1} = & C_{n,k,j+1} - \sum_{i=1}^{j+1} H_{k-1}^{(j+1-i)} (-1)^{j+1-i} c(n,i) \\  = & \sum_{i=1}^{j+2} (k-1) H_{k-2}^{(j+2-i)} (-1)^{j+2-i} c(n,i) - \sum_{i=1}^{j+1} H_{k-1}^{(j+1-i)} (-1)^{j+1-i} c(n,i) \\  = & (k-1)c(n,j+2) + \sum_{i=j+1} (-1)^{j+2-i} c(n,i) \left( (k-1) H_{k-2}^{(j+2-i)} + H_{k-1}^{(j+1-i)}\right) \\  = & (k-1) \left( c(n,j+2) + \sum_{i=j+1} (-1)^{j+2-i} c(n,i) \left( H_{k-2}^{(j+2-i)} + \frac{1}{k-1} H_{k-1}^{(j+1-i)}\right) \right) \\  = & (k-1) \left( c(n,j+2) + \sum_{i=j+1} (-1)^{j+2-i} c(n,i) H_{k-1}^{(j+2-i)} \right).  \end{align*}

But again by Theorem \ref{thm:main2}, $C_{n,k+1,j+1}$ is \[ \frac{k}{k-1}m_{n,k,j+1,j+1} = \sum_{i=j+2} (-1)^{j+2-i} k \cdot c(n,i) H_{k-1}^{(j+2-i)},\] as desired. $\square$

Loeb \cite{L} defines a generalization of the Stirling numbers $s(n,k)$ for arbitrary values of $n$, and examines these numbers when $n$ is a negative integer. Briefly, $s(n,k)$ can be defined for all integers $n$ and $k \geq 0$ using the recursive relation $s(n+1,k) = s(n,k-1) - ns(n,k)$ and the initial conditions $s(n,0) = \delta_{n0}$.

For our purposes, we will only use the following:

\begin{theorem} (\cite[Theorem 2]{L}) Using Loeb's \cite{L} definition of $s(-n,k)$, we have \[ (-1)^{k}n!s(-n,k) = H_{n}^{(k)}.\] \end{theorem}

Substituting this formula for $H_{m}^{(k)}$ into Theorem \ref{thm:lincomb} gives us an appealing formula for the coefficients $C_{n,k,j}$ and the $r$-Stirling numbers in terms of $c(n,k)$ and $s(-n,k)$.

\begin{corollary} For $j>0$ and $k>2$, we have
\[ C_{n,k,j} = \sum_{i=1}^{j+1} (k-1)!s(2-k,j+1-i) c(n,i) \] and \[ \rstir{n}{j+r}{r} = \sum_{i=1}^{j+1} s(1-r,j+1-i)c(n,i).\] \end{corollary}

Interestingly, this is not the only connection between $r$-Stirling numbers and iterated harmonic sums. We define another set of iterated harmonic sums as follows. Set \[ H_{n,j}^{1} = \sum_{1 \leq i_{1} < i_{2} < \cdots < i_{j} \leq n} \frac{1}{i_{1}i_{2} \cdots i_{j}}, \] and then for $k >1$, recursively define \[ H_{n,j}^{k} = \sum_{i=1}^{n} H_{i,j}^{k-1} .\] (Notice the slight change in notation to distinguish this notion from $H_{n}^{(j)}$.) Of course, $H_{n,j}^{k} = 0$ for any $n < j$.

First, we note that the generating function for the Stirling numbers of the first kind, \[ (x+1)(x+2) \cdots (x+n-1) = \sum_{j=0}^{n-1} c(n,j+1) x^{j}\] easily implies that \[ c(n,j+1) = \sum_{1 \leq i_{1} < i_{2} < \cdots < i_{n-j-1} \leq n-1} i_{1}i_{2} \cdots i_{n-j-1} = e_{n-j-1}(1,2, \ldots, n-1).\] (Here, $e_{k}$ is the usual elementary symmetric function.) We can re-write this formula for $j>0$ as \begin{equation} \label{eq:harm1} c(n,j+1) = (n-1)! \sum_{1 \leq i_{1} < i_{2} < \cdots < i_{j} \leq n-1} \frac{1}{i_{1}i_{2} \cdots i_{j}} = (n-1)! H_{n-1,j}^{1}.\end{equation}

Similarly, the $r$-Stirling numbers have generating function \[ (x+r)(x+r+1) \cdots (x+n-1) = \sum_{j=0}^{n-r} \rstir{n}{j+r}{r}x^{j}, \] so that \[ \rstir{n}{j+r}{r} = \sum_{r \leq i_{1} < i_{2} < \cdots < i_{n-j-r} \leq n-1} i_{1}i_{2} \cdots i_{n-j-r} = e_{n-j-r}(r,r+1,\ldots, n-1), \] which is equivalent to \[ = \frac{(n-1)!}{(r-1)!} \sum_{r \leq i_{1} < i_{2} < \cdots < i_{j} \leq n-1} \frac{1}{i_{1}i_{2} \cdots i_{j}}.\] A combinatorial interpretation is of this formula is described by Broder \cite[Theorem 7]{B}.

However, there is another way to generalize equation \eqref{eq:harm1}. We begin with another recurrence we can use to describe the $r$-Stirling numbers, which is due to Broder \cite[Lemma 11]{B}.

\begin{proposition} \label{prop:recur} For $r>1$ and $n \geq m \geq r$, the $r$-Stirling numbers satisfy \[ \rstir{n}{m}{r} = \sum_{i = 0}^{n-m}  \frac{(n-r)!}{(m+i-r)!} \rstir{m+i-1}{m-1}{r-1} .\] \end{proposition}

Proof. We count the permutations counted by $\rstir{n}{m}{r}$ according to the position of $n$. Consider a permutation $\sigma$ in $S_{n}$ with $m$ left-to-right maxima that include $n, n-1, \ldots, n-r+1$. Since $n$ is the $m$th largest left-to-right maximum, it must occur at position $m$ or later in $\sigma$. Now, for $i$ from $0$ to $n-m$, we notice that if $n$ is in position $m+i$, then the $n - i - m$ entries that occur after $n$ in $\sigma$ will not be left-to-right maxima, so they must be chosen from $1,2, \ldots, n-r$. There are then $\frac{(n-r)!}{(m+i-r)!}$ ways to choose and place these elements. We also note that exactly $m-1$ more left-to-right maxima must occur before $n$, including $n-1,n-2, \ldots, n-r+1$. This implies that $\textrm{red}(\sigma_{1}\sigma_{2} \cdots \sigma_{m+i-1})$ is a permutation in $S_{m+i-1}$ with exactly $m-1$ left-to-right maxima, including $m+i-1, m+i-2, \ldots, m+i-r+1$, of which there are exactly $\rstir{m+i-1}{m-1}{r-1}$. $\square$

Then another formula for the $r$-Stirling numbers is given by the following:

\begin{theorem} For $r \geq 1$, $m > r$ and $n \geq r$, \[ \rstir{n}{m}{r} = (n-r)! H_{n-r,m-r}^{r}.\] \end{theorem}

Proof. The result is clear if $r=1$, since then the formula coincides with that given in \eqref{eq:harm1}. We proceed by induction on $r$. Recall that \[(n-r)! H_{n-r,m-r}^{r} = (n-r)! \sum_{i=1}^{n-r} H_{i,m-r}^{r-1} = (n-r)! \sum_{i=m-r}^{n-r} H_{i,m-r}^{r-1} \] \[ = (n-r)! \sum_{i=0}^{n-m} H_{m+i-r, m-r}^{r-1} \]  We claim that, for $i$ from $0$ to $n-m$, the term $(n-r)! H_{m+i-r,m-r}^{r-1}$ counts the permutations $\sigma$ with $n$ in position $m+i$. By the proof of Proposition \ref{prop:recur}, there are $\frac{(n-r)!}{(m+i-r)!} \rstir{m+i-1}{m-1}{r-1}$ such permutations. However, inductively, we may assume that \[ \rstir{m+i-1}{m-1}{r-1} = (m+i-r)! H_{m+i-r,m-r}^{r-1}.\]

Then the total number of permutations we have just counted is \[  \frac{(n-r)!}{(m+i-r)!} \cdot (m+i-r)!H_{m+i-r,m-r}^{r-1}  = (n-r)! H_{m+i-r,m-r}^{r-1}.\] Then summing over all $i$ from $0$ to $n-m$ gives the desired result. $\square$

Although this particular description of this formula for the $r$-Stirling numbers seems to be new, the formula itself is not entirely new. We could trace back the nested sums $H_{n-r,m-r}^{r}$ to express $\rstir{n}{m}{r}$ as an integer linear combination of the original harmonic sums $H_{i,m-1}^{1} = \frac{1}{i!} \rstir{1+i}{m}{1}$. That formula appears to be a special case of a theorem of Broder \cite[Theorem 12, the case $p=1$]{B}.

Adamchik \cite{A} makes some other connections between the classical Stirling numbers and generalized harmonic numbers.

\section{Border Patterns}

Finally, we note that the $r$-Stirling numbers $\rstir{n}{r+2}{r}$ are also related to another permutation pattern. Kitaev and Liese \cite{KL}, introduce the \textit{border mesh pattern} $p$ in Figure \ref{fig:patternp}:
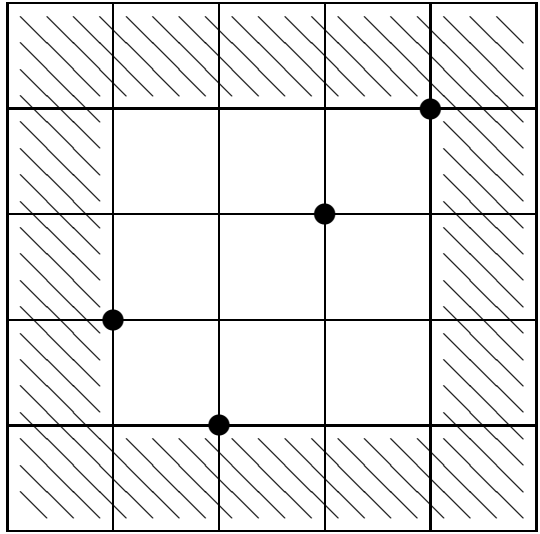
\begin{figure}[ht]
\begin{picture}(400,200)
\multiput(125,0)(40,0){6}{\line(0,1){200}}
\multiput(125,0)(0,40){6}{\line(1,0){200}}
\put(165,80){\circle*{8}}
\put(205,40){\circle*{8}}
\put(245,120){\circle*{8}}
\put(285,160){\circle*{8}}
\multiput(130,195)(10,0){14}{\line(1,-1){30}}
\put(270,195){\line(1,-1){50}}
\put(280,195){\line(1,-1){40}}
\put(290,195){\line(1,-1){30}}
\put(300,195){\line(1,-1){20}}
\put(310,195){\line(1,-1){10}}
\multiput(290,165)(0,-10){14}{\line(1,-1){30}}
\multiput(160,35)(10,0){13}{\line(1,-1){30}}
\multiput(130,185)(0,-10){13}{\line(1,-1){30}}
\put(130,55){\line(1,-1){50}}
\put(130,45){\line(1,-1){40}}
\put(130,35){\line(1,-1){30}}
\put(130,25){\line(1,-1){20}}
\put(130,15){\line(1,-1){10}}
\end{picture}
\caption{The pattern $p$ \label{fig:patternp}}
\end{figure}

This pattern is called a border pattern since the squares around the outside edge of the pattern are all marked ``$=0$''. We say that a permutation $\sigma = \sigma_{1} \cdots \sigma_{n}$ matches the pattern $p$ if $\sigma_{1} > 1$, $\sigma_{n} = n$, and there is an entry $a$ so that $\sigma_{1} < a < n$, and $a$ occurs after 1 in $\sigma$. (The first dot in the picture is $\sigma_{1}$, the second must be 1, and the fourth must be $n$. The third dot in the picture is $a$.) When we count occurrences of $p$, we will count every set of 4 entries of $\sigma$ that meet the description of the pattern. Since $\sigma_{1}$, 1, and $n$ must be involved in the pattern, this is equivalent to counting the number of possibilities for $a$ - i.e. the number of entries between 1 and $n$ that are larger than $\sigma_{1}$.

\begin{theorem} For any $k \geq 1$, the number of permutations in $S_{n}$ that match $p$ exactly $k$ times is exactly the number of permutations in $S_{n-1}$ that match $MMP^{k+1}$ exactly once, but almost match $MMP^{k+2}$. \end{theorem}

Proof. Again, we can give a bijective proof. First, we note that if $\sigma = \sigma_{1} \cdots \sigma_{n} \in S_{n}$ matches $p$ exactly $k$ times, then $\sigma = \sigma_{1} A 1 B n$, where $A$ and $B$ are strings, and $B$ contains exactly $k$ elements larger than $\sigma_{1}$. Then we form $\sigma' = \sigma_{1} B n A$, and let $\phi = (\sigma_{1}-1) B' (n-1) A'$, where $A'$ and $B'$ are obtained from $A$ and $B$, respectively, by subtracting 1 from each number in the string. Then $\phi \in S_{n-1}$ since $\sigma'$ contained $2, \ldots, n$.

Since the string $B$ contains exactly $k$ elements greater than  $\sigma_{1}$, $B'$ will contain exactly $k$ elements greater than $\sigma_{1}-1$. Hence $\phi$ will match $MMP^{k+1}$ but not $MMP^{k+2}$. This map will be a bijection since its inverse can be constructed easily. $\square$

\end{document}